\documentclass{article}

\usepackage{zharticle}
\usepackage{graphicx,amsmath,amsfonts,verbatim}

\newcommand{\ones}{{\bf 1{}}}  
\newcommand{\reals}{{\bf R{}}}  
\newcommand{\symms}{{\bf S{}}}  
\newcommand{\tr}{\mathop{\bf tr}}  
\newcommand{\norm}[1]{{\lVert #1 \rVert}}  

\newcommand{\eg}{{\rmfamily\itshape e.g.}}
\newcommand{\ie}{{\rmfamily\itshape i.e.}}

\newcommand{\cgr}{\texttt{cvxgenrust}}

\usepackage{lipsum}

\usepackage[colorlinks]{hyperref}

\usepackage{listings}
\usepackage{xcolor}

\definecolor{codegreen}{rgb}{0,0.6,0}
\definecolor{codegray}{rgb}{0.5,0.5,0.5}
\definecolor{codepurple}{rgb}{0.58,0,0.82}
\definecolor{backcolour}{rgb}{0.95,0.95,0.98}

\lstdefinestyle{zhstyle}{
    backgroundcolor=\color{backcolour},
    commentstyle=\color{codegreen},
    keywordstyle=\color{magenta},
    numberstyle=\tiny\color{codegray},
    stringstyle=\color{codepurple},
    basicstyle=\ttfamily\small,
    breakatwhitespace=false,
    breaklines=true,
    captionpos=b,
    keepspaces=true,
    numbers=left,
    numbersep=5pt,
    showspaces=false,
    showstringspaces=false,
    showtabs=false,
    tabsize=2,
}
\lstdefinelanguage{zhpython}{%
    sensitive=true,%
    morekeywords={and, as, assert, break, class, continue, def, del, elif,%
        else, except, exec, finally, for, from, global, if, import, in, is,%
        lambda, not, or, pass, print, raise, return, try, while, with, yield},%
    morecomment=[l]\#,%
    morestring=[s]{'''}{'''},%
    morestring=[s]{"""}{"""},%
    morestring=[b]',%
    morestring=[b]"%
}

\lstdefinelanguage{zhrust}{%
    sensitive=true,%
    morekeywords={as, async, await, break, const, continue, crate, dyn, else,%
        enum, extern, false, fn, for, if, impl, in, let, loop, match, mod, move,%
        mut, pub, ref, return, self, Self, static, struct, super, trait, true,%
        type, unsafe, use, where, while},%
    morecomment=[l]{//},%
    morecomment=[s]{/*}{*/},%
    morestring=[b]"%
}

\lstdefinestyle{Rust}{
    style=zhstyle,
    language=zhrust
}

\title{Generation of Custom Solvers in Rust\\for Convex Optimization}
\author[1,2]{Hao Zhu}
\author[1,2]{Joschka Boedecker}
\affil[1]{IMBIT//BrainLinks-BrainTools}
\affil[2]{Department of Computer Science, University of Freiburg}

\begin{document}
\maketitle

\begin{abstract}
We introduce \cgr, an open-source tool for generating custom Rust code that solves families of parameterized convex optimization problems modeled in CVXPY.
\cgr\ canonicalizes a problem family, extracts affine maps to Clarabel cone-program data, and generates a specialized Rust crate that updates parameters and calls Clarabel natively at runtime.
The generated solver can also be exposed to Python and registered as a custom CVXPY solver.
Our code generator supports a wide range of convex optimization problems up to semidefinite programs and exponential-cone problems.
Numerical experiments show reduced runtime relative to direct CVXPY solves and performance comparable to CVXPYgen on shared problem classes.
\end{abstract}

\newpage
\tableofcontents
\newpage

\section{Introduction}
We consider the class of \emph{parameterized convex optimization problems} of the form
\begin{equation}\label{prob:param_cp}
    \begin{array}{ll}
        \mbox{minimize} & f_0(x; \omega)\\
        \mbox{subject to} & f_i(x; \omega) \leq 0,\quad i = 1, \ldots, m\\
        & h_i(x; \omega) = 0,\quad i = 1, \ldots, p,
    \end{array}
\end{equation}
where $x \in \reals^n$ is the optimization variable~\cite{agrawal2019differentiable,schaller2022embedded}.
The \emph{family} of problems (\ref{prob:param_cp}) corresponding to different values of the parameter $\omega \in \reals^k$ is called a \emph{problem family}.
In this context, we call \emph{a} problem of the form (\ref{prob:param_cp}) associated with one specific value of $\omega$ a \emph{problem instance}.
The functions $f_i \colon \reals^n \times \reals^k \to \reals$ for $i = 0, \ldots, m$ are assumed to be convex and $h_i \colon \reals^n \times \reals^k \to \reals$ for $i = 1, \ldots, p$ are assumed to be affine, in the variable $x$ for each fixed $\omega$, so that each problem instance is a convex optimization problem.
Parameterized convex optimization problems are ubiquitous in many fields, including machine learning~\cite{mattingley2012cvxgen,agrawal2019differentiable,barratt2021least}, control~\cite{schaller2022embedded, ma2025learning}, and finance~\cite{boyd2024markowitz,luxenberg2024portfolio}, just to list a few.

In practice, solving a problem instance of the form (\ref{prob:param_cp}) typically involves two steps.
First, the problem instance is \emph{canonicalized} into a standard form, such as a cone program (see \S\ref{sec:quad_cone_prob}), that can be handled by a generic numerical solver.
Second, the canonicalized problem instance is passed to the numerical solver and solved via, \eg, an interior-point method~\cite{nesterov1994interior,boyd2004convex}.
Although the first problem canonicalization step is usually mathematically tedious, time-consuming, and error-prone, many \emph{domain specific languages} (DSLs) have been developed to automate this step, such as YALMIP~\cite{lofberg2004yalmip} and CVX~\cite{grant2014cvx} for MATLAB, Convex.jl~\cite{udell2014convex} for Julia, CVXPY~\cite{diamond2016cvxpy} for Python, and CVXR~\cite{fu2020cvxr} for R.
These DSLs allow users to express optimization problems in a natural mathematical syntax, and automatically canonicalize them into a standard form that can be solved by generic solvers such as OSQP~\cite{stellato2020osqp}, SCS~\cite{odonoghue2016conic}, ECOS~\cite{domahidi2013ecos}, and Clarabel~\cite{goulart2024clarabel}.

Most of the DSLs for convex optimization are \emph{parser-solvers}~\cite{schaller2022embedded}, which means that they parse an instance of the problem (\ref{prob:param_cp}) and call a generic solver at runtime, for each time the parameter $\omega$ changes.
This approach is flexible and easy to use, but it can be inefficient when the problem parameter $\omega$ changes frequently, as is often the case in real-time embedded applications. In such settings, the parsing and canonicalization steps may become computationally expensive, especially for large-scale problems.
Here we address this issue by designing a \emph{code generator}, which parses a given \emph{problem family} and generates custom solvers for it.
The generated source code is then compiled into binary libraries that are directly executable and can solve problem instances when the parameter values change, without the need for parsing and canonicalization at runtime.

\subsection{Prior and related work}
The idea of code generators for embedded convex optimization has been explored over the last several decades.
One of the earliest works in this direction is CVXGEN~\cite{mattingley2012cvxgen}, which generates custom C code for solving small linear programs (LPs) and convex quadratic programs (QPs) in embedded applications.
This software has been used to generate flight code for high-speed onboard convex optimization for precision landing of space rockets by SpaceX~\cite{blackmore2016autonomous}.
More recently, the CVXPY-based C code generator CVXPYgen~\cite{schaller2022embedded} extends the code generation approach to a much broader class of convex optimization problems, up to second-order cone programs (SOCPs).
Custom solvers generated by CVXPYgen are generally faster than those from CVXGEN, and have smaller binary sizes (see~\cite{schaller2022embedded} for detailed comparisons).
There are several other code generators, which, however, only support solver generation for a specific problem family, \eg, optimal control problems, including FORCESPRO~\cite{embotech2026forcespro}, FORCES~NLP~\cite{zanelli2020forces}, \texttt{acados}~\cite{verschueren2022acados}, and SCvxPyGen~\cite{berrah2024scvxpygen}.

There are several limitations of the existing code generators for convex optimization.
Firstly, existing code generators target C, mainly because C offers native performance, small runtime overhead, and broad support for embedded deployment.
In this work, we target Rust in order to retain these systems-level advantages while improving the safety and maintainability of the generated solver code.
Rust compiles to efficient native code without a garbage collector, and its ownership model helps prevent memory errors, unintended aliasing, and data races at compile time.
These properties make Rust a suitable target for generated solvers used in real-time and embedded applications.
Secondly, as far as we know, there is no existing code generator that supports solver generation for semidefinite programs (SDPs) or exponential cone programs, which have become increasingly useful in many applications (see \S\ref{sec:eg_sdp} and \S\ref{sec:eg_ml} for some examples).

\subsection{Contribution}
This work presents \cgr, a code generator for parameterized convex optimization problems that generates custom solvers in Rust. Given a problem family specified in CVXPY, \cgr\ canonicalizes the family once into a standard conic form and generates a Rust crate specialized to that family. At runtime, users only update the problem parameters and call the compiled solver, avoiding repeated parsing and canonicalization.

The generated solver uses Clarabel as its backend conic solver. Since Clarabel is implemented natively in Rust, \cgr\ integrates directly with the Rust ecosystem without foreign-function interfaces or language bindings between the generated solver and the numerical backend. This enables generation of native Rust solver modules that can be used from Rust projects, compiled into binary libraries, or exposed to other languages.

Compared with existing general-purpose code generators for convex optimization problems, \cgr\ targets Rust and supports a broader class of cone programs through Clarabel, including LPs, QPs, SOCPs, SDPs, and exponential cone programs. To the best of our knowledge, \cgr\ is the first code generator for parameterized convex optimization that combines native Rust solver generation with support for semidefinite and exponential cones.

In addition to the Rust module, \cgr\ generates Python bindings that allow the generated solver to be registered as a custom solver in CVXPY. Users can therefore keep the usual CVXPY modeling workflow while reducing runtime overhead for repeated solves of the same problem family.

We evaluate \cgr\ on problem families from several application domains. The numerical results show that generated \cgr\ solvers reduce runtime compared with standard CVXPY workflows and achieve solve times comparable to CVXPYgen on problem classes supported by both tools.

Last but not least, \cgr\ is fully open-source and available at
\begin{quote}
    \url{https://github.com/dxogrp/cvxgenrust}.
\end{quote}

\subsection{Outline}
The rest of this paper is organized as follows.
In \S\ref{sec:cgr} we provide a high-level overview of the basic ideas behind \cgr, followed by a simple example of how these ideas work and how to use \cgr\ for code generation.
In \S\ref{sec:num_results} we present numerical results on the performance of \cgr\ on a variety of problem families.
(The detailed problem formulations of the compared problem families can be found in the appendix \S\ref{sec:bench_prob}.)
We conclude the paper in \S\ref{sec:conclusion} with a summary and discussion of future work.

\section{\cgr}\label{sec:cgr}
\subsection{Canonical quadratic cone problems}\label{sec:quad_cone_prob}
CVXPY allows users to describe convex optimization problems in a high-level mathematical syntax, but numerical solvers operate on canonical forms.
For the Clarabel backend used by \cgr, each CVXPY problem instance is canonicalized into a quadratic cone program
\begin{equation}\label{prob:cone_qp}
    \begin{array}{ll}
        \mbox{minimize} & (1/2)\tilde{x}^T P\tilde{x}
            + q^T \tilde{x}\\
        \mbox{subject to} & A\tilde{x} + s = b\\
        & s \in \mathcal{K},
    \end{array}
\end{equation}
where $\tilde{x}$ is the canonical optimization variable, $s$ is the slack variable, and $\mathcal{K}$ is a product cone.
The data $P$, $q$, $A$, and $b$ are the canonical numerical data passed to Clarabel, while $\mathcal{K}$ is determined automatically by the atoms and constraints in the original CVXPY model.
Here the supported cone components include zero, nonnegative, second-order, positive semidefinite, exponential, and power cones.

\subsection{Affine-solve-affine workflow}
We could similarly define the data, or parameters, of the canonicalized cone program (\ref{prob:cone_qp}) as $\tilde{\omega} = \{P, q, A, b\}$.
The important point for code generation is that, for \emph{disciplined parametrized programming} (DPP) compliant problems~\cite{agrawal2019differentiable}, the canonical problem parameters $\tilde{\omega}$ depend affinely on the parameters $\omega$ of the original problem (\ref{prob:param_cp}), and the original problem solution $x^\star$ can be recovered from the canonical solution $\tilde{x}^\star$ via an affine map.
In other words, we have
\begin{equation*}
    \tilde{\omega} = F \left[\begin{array}{c}
        \omega \\ 1
    \end{array}\right],\qquad x^\star = G \left[\begin{array}{c}
        \tilde{x}^\star \\ 1
    \end{array}\right],
\end{equation*}
where $F$ and $G$ are sparse matrices depending only on the structure of the problem (\ref{prob:param_cp}) and can be extracted from the CVXPY-Clarabel canonicalization process.

\begin{figure}[p]
    \centering
    \includegraphics[width=0.65\textwidth]{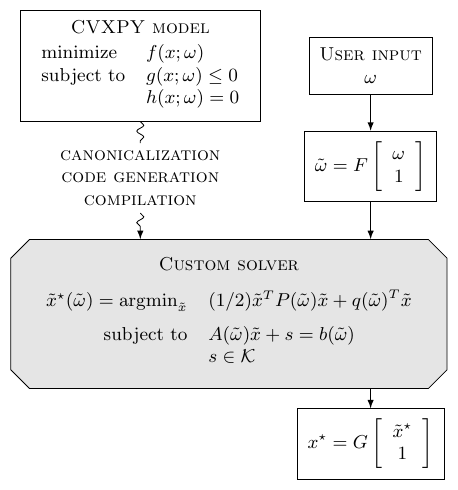}
    \caption{
        Workflow of \cgr.
        The left part shows the offline code generation process and the right part shows the runtime solving process.
    }\label{fig:workflow}
\end{figure}

This \emph{affine-solve-affine} workflow~\cite{agrawal2019differentiable} separates offline code generation from runtime solving.
Offline, \cgr\ extracts the canonicalization maps; at runtime, it updates the parameters, reconstructs the cone program, calls Clarabel, and maps the solution back.
Figure~\ref{fig:workflow} summarizes this workflow.

\subsection{A simple example}
As a minimal example, consider the nonnegative least-squares problem
\begin{equation}\label{prob:nnls}
    \begin{array}{ll}
        \mbox{minimize} & \norm{Cx - d}_2^2\\
        \mbox{subject to} & x \succeq 0,
    \end{array}
\end{equation}
with parameters $C \in \reals^{m \times n}$ and $d \in \reals^m$ ($\succeq$ denotes componentwise inequality).
Introducing the residual variable $t = Cx - d$, we can rewrite (\ref{prob:nnls}) in the form of (\ref{prob:cone_qp}) with canonical variable $\tilde{x} = \left[\begin{array}{c}
    x \\ t
\end{array}\right] \in \reals^{n + m}$ and slack variable $s = \left[\begin{array}{c}
    s_1 \\ s_2
\end{array}\right] \in \reals^{m + n}$.
The corresponding problem data are
\begin{equation*}
    P = \left[\begin{array}{cc}
        0 & 0\\
        0 & 2I
    \end{array}\right],\qquad
    q = 0,\qquad
    A = \left[\begin{array}{cc}
        C & -I\\
        -I & 0
    \end{array}\right],\qquad
    b = \left[\begin{array}{c}
        d\\
        0
    \end{array}\right],
\end{equation*}
with cone $\mathcal{K} = {\{0\}}^m \times \reals_+^n$ (where the first block enforces $Cx - t = d$, and the second block enforces $x \succeq 0$).
In this example, the canonical problem parameters $P$, $q$, $A$, and $b$ are indeed affine in the original problem parameters $C$ and $d$.

The following CVXPY code generates a Rust solver crate for this problem family.
\begin{lstlisting}[language=zhpython]
import cvxpy as cp
import cvxgenrust as cgr

# Define the problem family
x = cp.Variable(n, name="x")
C = cp.Parameter((m, n), name="C")
d = cp.Parameter(m, name="d")

problem = cp.Problem(
    cp.Minimize(cp.sum_squares(C @ x - d)),
    [x >= 0]
)

# Generate the custom solver code
cgr.generate_code(
    problem,
    module_name="nonneg_ls"
)
\end{lstlisting}
The names assigned to CVXPY parameters and variables are used to generate stable parameter setters and solution extractors.
After code generation, the resulting Cargo project can be built and used from Rust.
For example, a Rust application can create a generated problem object, set the current values of $C$ and $d$, solve the problem, and extract the named variable $x$:
\begin{lstlisting}[style=Rust]
use nonneg_ls::CGRProblem;

// Create a new problem instance
let mut problem = CGRProblem::new();

// Update the problem parameters
problem.set_c(&c_values)?;
problem.set_d(&d_values)?;

// Solve and extract the solution
let solution = problem.solve()?;
let x = problem.extract_x(&solution.x)?;
\end{lstlisting}
To register and use the generated solver in CVXPY, the generated Python bindings can be used as follows:
\begin{lstlisting}[language=zhpython]
from nonneg_ls_wrapper.cgr_solver import cgr_solve

# Register the generated solver in CVXPY
problem.register_solve("CGR", cgr_solve)

# Update the problem parameters
C.value = C_value
d.value = d_value

# Solve with the generated solver
problem.solve(
    method="CGR",
    updated_params=["C", "d"]
)
\end{lstlisting}
\section{Numerical results}\label{sec:num_results}
We benchmarked \cgr\ on the twelve problem families listed in \S\ref{sec:bench_prob}, with five problem sizes per family.
Each case used three warm-up solves and 1000 timed samples; all methods used Clarabel with the same tolerances where applicable.
Experiments ran on Linux with an AMD Ryzen~9 7950X CPU.

We compare direct CVXPY solves, the generated \cgr\ solver called through its Python wrapper, the generated \cgr\ Rust crate called directly, and the corresponding Python and C interfaces generated by CVXPYgen.
Note that CVXPYgen is not available for the parametric quadratic-form, SDP, and logistic-regression examples in this benchmark suite.

\begin{figure}
    \centering
    \includegraphics[width=\textwidth]{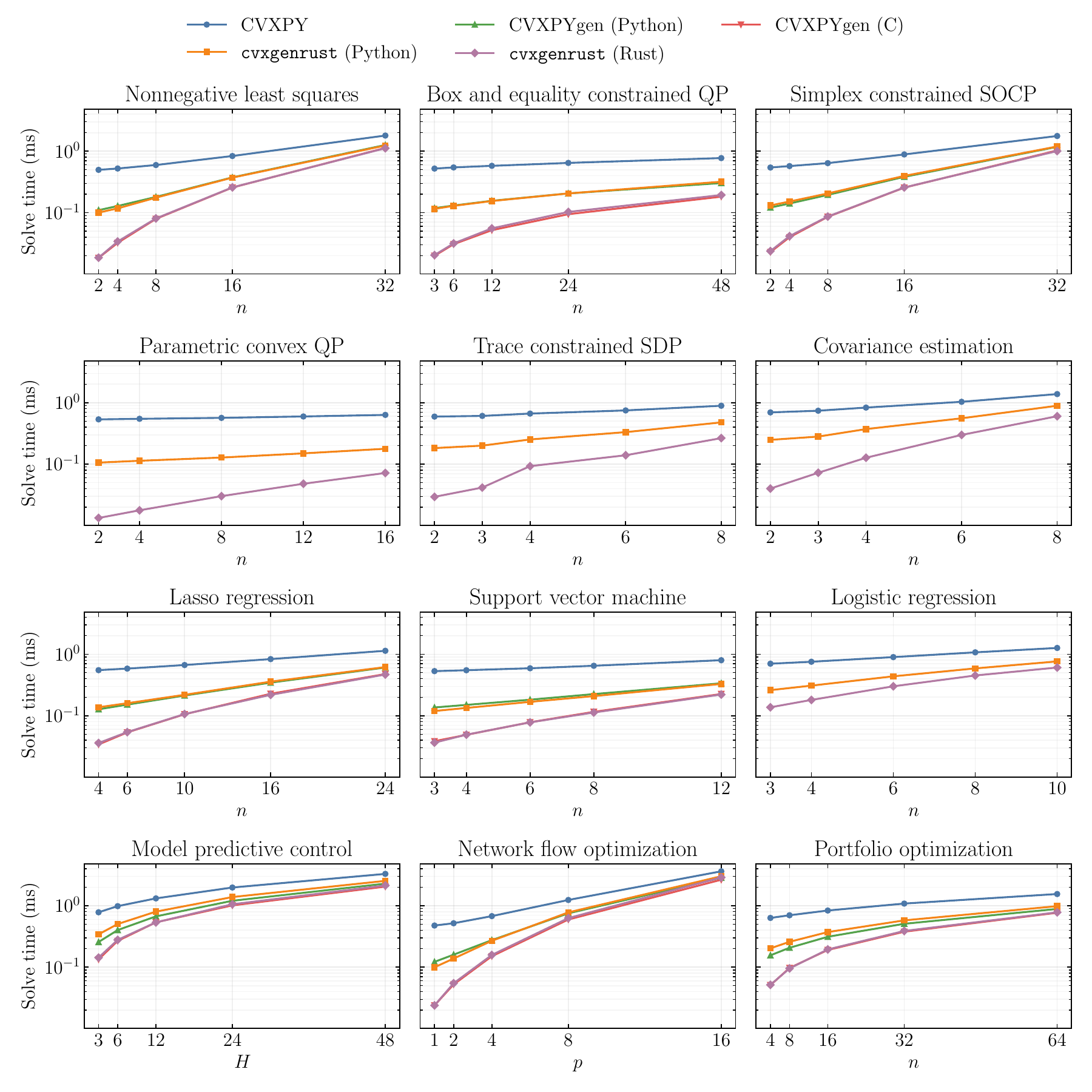}
    \caption{
        Comparison of solve times for CVXPY, CVXPYgen, and \cgr\ (shown in different colors).
        Each point represents the median solve time for one problem size over the timed samples.
    }\label{fig:solve_times}
\end{figure}

Figure~\ref{fig:solve_times} shows that both generated \cgr\ interfaces reduce runtime relative to direct CVXPY solves, with median speedups of $2.66\times$ through Python and $6.27\times$ through direct Rust calls over the 60 problem-size cases.
This improvement comes from avoiding repeated CVXPY parsing and canonicalization after code generation.
On the 40 problem-size cases supported by CVXPYgen, the \cgr\ Python and Rust interfaces are close to the corresponding CVXPYgen Python and C interfaces.

\section{Conclusion}\label{sec:conclusion}
We presented \cgr, a code generator that turns DPP-compliant CVXPY problem families into specialized Rust solver crates.
By extracting affine canonicalization maps offline and calling Clarabel natively at runtime, \cgr\ avoids repeated parsing and canonicalization while supporting a broad class of cone programs, including LPs, QPs, SOCPs, SDPs, and exponential-cone problems.
Numerical experiments show reduced runtime compared with direct CVXPY solves and performance comparable to CVXPYgen on shared problem classes.

Future work involves extending \cgr\ to support differentiating through the generated solvers (\eg, as in \cite{healey2026differentiating}), which would enable end-to-end learning of problem parameters in applications such as machine learning and control.
Moreover, we could extend \cgr\ to support generating explicit solver code for specific problem families, \eg, a parameterized QP~\cite{schaller2025automatic}, which can be solved more efficiently than general cone programs.

\section{Appendix: Benchmark problem families}\label{sec:bench_prob}
In this section we present the problem families used in the numerical experiments.
For the detailed problem dimensions and parameter generation procedures, please refer to the code in our repository
\begin{quote}
    \url{https://github.com/dxogrp/cgr-benchmarks}.
\end{quote}

\subsection{Basic examples}
\paragraph{Nonnegative least squares.}
We consider the problem of finding a nonnegative vector $x \in \reals^n$ that minimizes the least squares error with respect to a given matrix $A \in \reals^{m \times n}$ and a vector $b \in \reals^m$:
\begin{equation*}
    \begin{array}{ll}
        \mbox{minimize} & \norm{Ax - b}_2^2 \\
        \mbox{subject to} & x \succeq 0,
    \end{array}
\end{equation*}
where $\succeq$ represents elementwise inequality.
This problem is parameterized by the matrix $A$ and the vector $b$.

\paragraph{Box and equality constrained QP.}
We consider the following QP with box and equality constraints:
\begin{equation*}
    \begin{array}{ll}
        \mbox{minimize} & (1/2)\norm{x - x_0}_2^2 + q^T x \\
        \mbox{subject to} & Ax = b\\
        & 0 \preceq x \preceq 1
    \end{array}
\end{equation*}
where $x \in \reals^n$ is the optimization variable, and $x_0 \in \reals^n$, $q \in \reals^n$, $A \in \reals^{m \times n}$, and $b \in \reals^m$ are parameters of the problem.

\paragraph{Simplex constrained SOCP.}
We consider the following SOCP with the probabilistic simplex constraint:
\begin{equation*}
    \begin{array}{ll}
        \mbox{minimize} & \norm{Ax - b}_2^2 + 0.1 \norm{x}_2^2 \\
        \mbox{subject to} & \norm{x}_2 \leq \rho\\
        & x \succeq 0,\quad \ones^T x = 1,
    \end{array}
\end{equation*}
where $x \in \reals^n$ is the variable, and $A \in \reals^{m \times n}$, $b \in \reals^m$, and $\rho \in \reals_+$ are the parameters.

\paragraph{Parametric convex QP.}
We consider the following convex QP with the quadratic coefficient matrix $P \in \symms_+^n$ and a vector $q \in \reals^n$ as parameters:
\begin{equation*}
    \begin{array}{ll}
        \mbox{minimize} & x^T P x + q^T x \\
        \mbox{subject to} & x \succeq 0,\quad \ones^T x \leq 1,
    \end{array}
\end{equation*}
where $x \in \reals^n$ is the variable.
Note that \cgr\ allows modeling of this problem family, but it is not yet supported by CVXPYgen~\cite{schaller2022embedded}.

\subsection{Semidefinite programming examples}\label{sec:eg_sdp}
\paragraph{Trace constrained SDP.}
We consider the following SDP with equality constraints:
\begin{equation*}
    \begin{array}{ll}
        \mbox{minimize} & \tr(C^T X)\\
        \mbox{subject to} & X \succeq 0,\quad \tr X = 1,
    \end{array}
\end{equation*}
where $X \in \symms^n$ is the variable, and $C \in \symms^n$ is the parameter.
Here the inequality $\succeq$ denotes matrix inequality, \ie, $X$ is constrained to be positive semidefinite.

\paragraph{Covariance estimation.}
We consider the problem of recovering a covariance matrix $X \in \symms^n_+$ from a given noisy observation $S \in \symms^n$:
\begin{equation*}
    \begin{array}{ll}
        \mbox{minimize} & \norm{X - S}^2_F + 0.15\sum_{i \neq j} {(X_{ij} - S_{ij})}^2 \\
        \mbox{subject to} & X \succeq 0\\
        & X_{ii} = S_{ii},\quad i = 1, \ldots, n,
    \end{array}
\end{equation*}
where $\norm{\,\cdot\,}_F$ denotes the Frobenius norm, $X_{ij}$ and $S_{ij}$ are the $(i,j)$th entries of $X$ and $S$, respectively, and the problem is parameterized by the noisy observation $S$.

\subsection{Machine learning examples}\label{sec:eg_ml}
\paragraph{Lasso regression.}
The lasso regression problem is given by
\begin{equation*}
    \begin{array}{ll}
        \mbox{minimize} & (1/2)\norm{Ax - b}_2^2 + \lambda \norm{x}_1,
    \end{array}
\end{equation*}
where $x \in \reals^n$ is the variable, and $A \in \reals^{m \times n}$, $b \in \reals^m$, and $\lambda \in \reals_+$ are the parameters.

\paragraph{Support vector machine.}
We consider the soft-margin support vector machine problem:
\begin{equation*}
    \begin{array}{ll}
        \mbox{minimize} & (1/2)\norm{w}_2^2 + \sum_{i = 1}^{m} {(1 - y_i (X_i w + \beta))}_+,
    \end{array}
\end{equation*}
where $w \in \reals^n$ and $\beta \in \reals$ are the variables, and $X \in \reals^{m \times n}$ and $y \in {\{-1, 1\}}^m$ are the parameters.
The function ${(u)}_+$ denotes the positive part of $u \in \reals$.

\paragraph{Logistic regression.}
The logistic regression problem is given by
\begin{equation*}
    \begin{array}{ll}
        \mbox{minimize} & (1/m) \sum_{i = 1}^{m} \log(1 + \exp(-y_i (X_i w + \beta))) + (\lambda/2) \norm{w}_2^2,
    \end{array}
\end{equation*}
where $w \in \reals^n$ and $\beta \in \reals$ are the variables, and $X \in \reals^{m \times n}$, $y \in {\{-1, 1\}}^m$, and $\lambda \in \reals_+$ are the parameters.
Note that code generation for this problem family requires canonicalizing it into a conic form involving the exponential cone, which is currently not compatible with CVXPYgen.

\subsection{Control and finance examples}
\paragraph{Model predictive control.}
We consider the following finite-horizon linear-quadratic tracking problem.
Let $x_t \in \reals^3$ be the system state at time $t$ for $t = 0, \ldots, H$, and let $u_t \in \reals^2$ be the control input at time $t$ for $t = 0, \ldots, H-1$.
Given an initial state $\bar{x}_0 \in \reals^3$ and a reference trajectory $r_t \in \reals^3$ for $t = 0, \ldots, H$, we consider the following optimization problem:
\begin{equation*}
    \begin{array}{ll}
        \mbox{minimize} & \sum_{t = 0}^{H} \norm{Q^{1/2}(x_t - r_t)}_2^2 + \sum_{t = 0}^{H-1} \norm{R^{1/2} u_t}_2^2 + 0.03\norm{x_H}_2^2 \\
        \mbox{subject to} & x_0 = \bar{x}_0,\quad x_{t+1} = A x_t + B u_t,\quad t = 0, \ldots, H-1\\
        & |u_{t+1} - u_t| \leq 0.35,\quad t = 0, \ldots, H-2\\
        & |x_{1, t}| \leq 2.5,\quad |x_{2, t}| \leq 2.5,\quad t = 0, \ldots, H\\
        & |u_t| \leq 1,\quad t = 0, \ldots, H-1,
    \end{array}
\end{equation*}
where $x_t$ and $u_t$ are the optimization variables.
In this example we assume that the system dynamics $A \in \reals^{3 \times 3}$ and $B \in \reals^{3 \times 2}$, as well as the quadratic cost matrices $Q \in \symms_+^3$ and $R \in \symms_+^2$, are constants, so the problem is only parameterized by the initial state $\bar{x}_0$ and the reference trajectory $r_t$ for $t = 0, \ldots, H$.

\paragraph{Network flow optimization.}
We consider a capacitated multi-commodity flow problem on a directed graph with $n$ vertices and $m$ edges.
Let $p$ denote the number of commodities, and for each commodity $i = 1, \ldots, p$, let $f_i \in \reals^m$ represent the flow of the $i$th commodity on each directed edge, so $f = {\left[\begin{array}{ccc}
    f_1 & \cdots & f_p
\end{array}\right]}^T \in \reals^{p \times m}$ is the commodity-edge flow matrix.
The network flow optimization problem is given by
\begin{equation*}
    \begin{array}{ll}
        \mbox{minimize} & \sum_{i = 1}^{p} c_i^T f_i\\
        \mbox{subject to} & f_i \succeq 0,\quad Bf_i = d_i b_i,\quad i = 1, \ldots, p\\
        & \sum_{i = 1}^{p} f_{ij} \leq u_j,\quad j = 1, \ldots, m,
    \end{array}
\end{equation*}
where $c_i \in \reals^m$ is the edge cost vector for the $i$th commodity, $B \in \reals^{n \times m}$ is the vertex-edge incidence matrix of the graph, \ie, $B_{ij} = 1$ if edge $j$ enters vertex $i$, $B_{ij} = -1$ if edge $j$ leaves vertex $i$, and $B_{ij} = 0$ otherwise.
The vector $d \in \reals_+^p$ records the demand for all commodities, and the vector $b_i \in \reals^n$ encodes the source and sink of this commodity, which has value $1$ at the source vertex, value $-1$ at the sink vertex, and value $0$ elsewhere.
Therefore the conservation constraint $Bf_i = d_i b_i$ ensures that the flow of the $i$th commodity satisfies the demand $d_i$ from the source to the sink.
The vector $u \in \reals_+^m$ is the edge capacity vector.
This problem has the variable $f \in \reals^{p \times m}$ and is parameterized by $c_1, \ldots, c_p \in \reals^m$, $d \in \reals_+^p$, and $u \in \reals_+^m$.
The matrix $B \in \reals^{n \times m}$ and the vectors $b_1, \ldots, b_p \in \reals^n$ are assumed to be constants.

\paragraph{Portfolio optimization.}
We consider the long-only portfolio optimization problem with $n$ assets, given by
\begin{equation*}
    \begin{array}{ll}
        \mbox{minimize} & \gamma\left(\norm{F^T w}_2^2 + \norm{D^{1/2} w}_2^2\right) - \mu^T w + 0.02\norm{w - w^{\rm prev}}_1\\
        \mbox{subject to} & 0 \preceq w \preceq w^{\rm max},\quad \ones^T w = 1\\
        & \norm{w - w^{\rm prev}}_1 \leq 0.75,\quad Sw \preceq 0.65\ones,\quad \mu^T w \geq r^{\rm tar},
    \end{array}
\end{equation*}
where the wealth allocation vector $w \in \reals^n$ is the variable.
The matrix $F \in \reals^{n \times p}$ is the factor loading matrix with $p$ factors, and $D \in \symms_+^n$ is the diagonal idiosyncratic risk matrix, so $\gamma \geq 0$ is the risk aversion parameter that controls the trade-off between risk and return in the objective.
The vector $\mu \in \reals^n$ is the expected return vector, and $w^{\rm prev} \in \reals^n$ is the previous wealth allocation.
The matrix $S \in \reals^{q \times n}$ encodes the sector information of the assets, where $q$ is the number of sectors and $S_{ij} = 1$ if asset $j$ belongs to sector $i$ and $S_{ij} = 0$ otherwise.
The vector $w^{\rm max} \in \reals_+^n$ is the maximum allocation vector, and $r^{\rm tar} \in \reals$ is the target return.
We assume that the problem is parameterized by $\gamma$, $\mu$, $w^{\rm prev}$, and $r^{\rm tar}$, while the matrices $F$, $D$, and $S$, as well as the vector $w^{\rm max}$, are constants.

\newpage
\section*{Acknowledgments}
This work has been funded as part of BrainLinks-BrainTools, which is funded by the Federal Ministry of Economics, Science and Arts of Baden-W\"urttemberg within the sustainability program for projects of the Excellence Initiative II, and CRC/TRR 384 ``IN-CODE''.

\newpage
\bibliography{refs}

\end{document}